\documentclass[reqno]{amsart}

\pdfoutput=1

\usepackage{amsfonts,amsmath,amsthm,amsxtra,amssymb,mathrsfs,dsfont}

\usepackage{tikz-cd}

\usepackage{graphicx}  

\usepackage{subfigure}

\usepackage{wrapfig}

\usepackage{fancybox}

\usepackage[colorlinks = true, citecolor = blue, urlcolor = blue]{hyperref}

\theoremstyle{definition}

\def\gor#1{\widetilde{#1}}

\def\h(#1,#2){\mbox{Hom}\left(#1,#2\right)}

\def\t(#1,#2){\mbox{Tor}\left(#1,#2\right)}
\def\e(#1,#2){\mbox{Ext}\left(#1,#2\right)}

\def\F{\mathbb{F}}

\def\N{\mathbb{N}}

\def\R{\mathbb{R}}

\def\X{\mathbb{X}}

\def\Z{\mathbb{Z}}

\def\uu{\mathbf{u}}
\def\vv{\mathbf{v}}

\def\dgm{\mathsf{bcd}}
\def\Dgm{\mathsf{Bcd}}

\def\<{\langle}
\def\>{\rangle}

\title{A Brief History of Persistence}
\author[Jose Perea]{Jose A. Perea }
\address{
\shortstack[l]{
Department of Computational Mathematics, Science \& Engineering \\
Department of Mathematics, \\
Michigan State University \\
East Lansing, MI, USA.}}
\email{joperea@msu.edu}

\begin{document}
\makeatletter{\renewcommand*{\@makefnmark}{}
\footnotetext{This work was partially supported by the NSF (DMS-1622301) and  DARPA (HR0011-16-2-003)}\makeatother}
\makeatletter{\renewcommand*{\@makefnmark}{}
\makeatother}

\maketitle

\begin{abstract}
Persistent homology is currently  one of the more widely known tools from computational topology and topological data analysis. We present in this note a brief survey on the evolution of the subject.
The goal is to highlight the main ideas, starting from the subject's computational inception more than 20 years ago,
to the more modern categorical and representation-theoretic point of view.
\end{abstract}

\section{The Early Years: Computer Science Meets Geometry, Algebra and Topology}
Our current view of persistent homology can be traced back to
work  of Patrizio Frosini (1992) on size functions \cite{frosini1992measuring},
and of Vanessa Robins  (1999) \cite{robins1999towards}
on using experimental data
to infer  the topology of attractors
in dynamical systems.
Both  approaches rely on singular homology  as a shape descriptor, which leads to what is known today as
the \emph{``homology inference problem''}:
 Given a finite set $X$ (the data) sampled from/around a topological space $\X$ (e.g., the attractor),
how can one infer the homology of $\X$ from $X$ with high confidence?
See for instance \cite{niyogi2008finding} for  the case when
$\X$ is a compact Riemannian submanifold of Euclidean space, and
$X \subset \X$ is sampled according to the intrinsic uniform distribution.
From here on out it will be useful to think of $X$ and $\X$ as  subspaces of a  bounded
metric space $(\mathbb{M}, \rho)$. In this case, one can formalize  the statement ``$X$ approximates $\X$''
 by saying that if
$Z\subset \mathbb{M}$, $\epsilon \geq 0$, and  $Z^{(\epsilon)} := \{x\in \mathbb{M} : \rho(x, Z) \leq \epsilon\}$, then the Hausdorff distance
\[
d_H (X,\X) := \inf \left\{\epsilon > 0 : X \subset \X^{(\epsilon)} \mbox{ and }
\X \subset X^{(\epsilon)}\right\}
\]
is small.
The goal is then to approximate the topology of $\X$ by that of $X^{(\epsilon)}$.
Below in Figure \ref{fig:OffSetFiltration} we illustrate the evolution of $X^{(\epsilon)}$ as $\epsilon$ increases.
\begin{figure}[htb!]
  \centering
  \includegraphics[width= \textwidth]{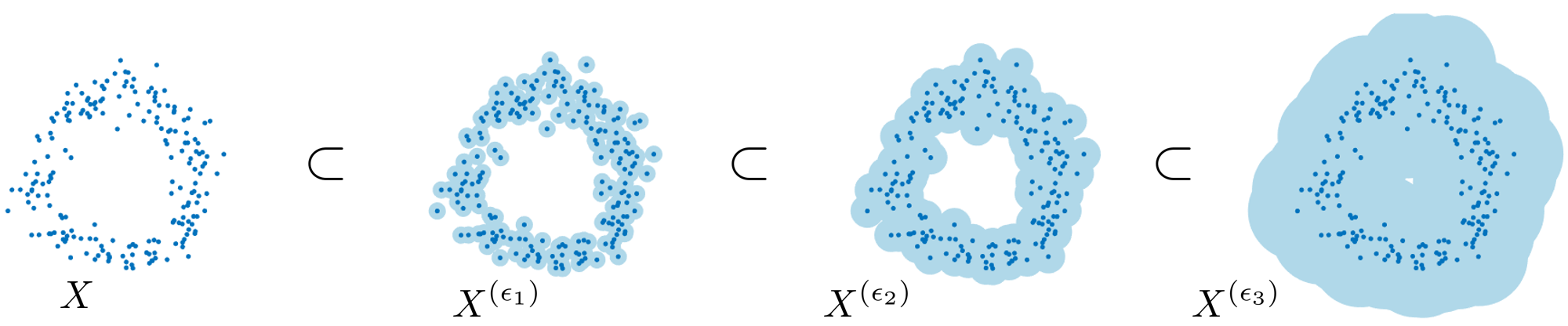}
  \caption{Some examples of $X^{(\epsilon)}$ for $X\subset \R^2$ sampled around the unit circle, and $\epsilon$ values $0 < \epsilon_1 < \epsilon_2< \epsilon_3$. }
  \label{fig:OffSetFiltration}
\end{figure}

In order to capture the multiscale nature of $\mathcal{X} = \{X^{(\epsilon)}\}_{\epsilon }$,
and deal with the instability of topological features in $X^{(\epsilon)}$ as $\epsilon$ changes,
Frosini and Robins introduced (independently) the idea of \emph{homological persistence}:
for   $\epsilon, \delta \geq 0$ let
\[
\iota^{\epsilon,\delta}: X^{(\epsilon)} \hookrightarrow X^{(\epsilon + \delta)}
\]
be the inclusion,
and  consider the  induced
linear map in homology  with coefficients in a field $\F$
\[
\iota^{\epsilon,\delta}_* :
H_n\left(X^{(\epsilon)} ; \F\right)
\longrightarrow
H_n\left(X^{(\epsilon + \delta)}; \F \right)
\]
The image of $\iota^{\epsilon,\delta}_*$ is the $\delta$-persistent $n$-th homology group of
the filtered space
$\mathcal{X}$ at $\epsilon$,
denoted $H_n^{\epsilon,\delta}(\mathcal{X}; \F)$;
and $\mathsf{rank}\left(\iota_*^{\epsilon,\delta}\right)$
is the persistent
Betti number
$\beta^{\epsilon, \delta}_n(\mathcal{X}; \F) $.

The design of algorithms to efficiently compute/approximate these integers
is of course predicated on first
replacing the spaces $X^{(\epsilon)}$ by   finite, combinatorial  models of their topology.
Fortunately there is a vast literature  on how to do this.
Take for instance the Vietoris-Rips  complex,
first introduced by Leopold Vietoris in the nineteen-twenties
in an attempt to define a homology theory for general metric spaces \cite{vietoris1927hoheren}.
It is defined, for $Z\subset \mathbb{M}$  and $\epsilon \geq 0$,
as the abstract simplicial complex
\[
R_\epsilon(Z) := \big\{ \{z_0,\ldots, z_k\} \subset Z: \rho(z_i, z_j) \leq \epsilon \; \mbox{ for all }\; 0 \leq i,j \leq k\big\}
\]
Below in Figure \ref{fig:RipsFiltration} we show an example of how $R_\epsilon(Z)$ evolves as $\epsilon$
increases, for $Z\subset \R^2$ sampled around the unit  circle,
and for $\epsilon$ values $0< \epsilon_1  < \epsilon_2 < \epsilon_3$.
\begin{figure}[htb!]
  \centering
  \includegraphics[width= \textwidth]{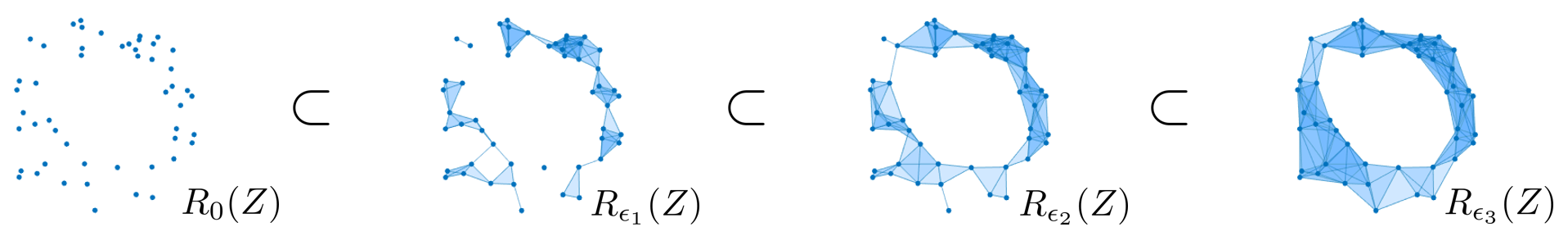}
  \caption{Some examples of the Rips complex, for points sampled around the unit circle in $\R^2$.}
  \label{fig:RipsFiltration}
\end{figure}

Notice that  $R_\epsilon(Z) \subset R_{\epsilon + \delta} (Z)$ whenever $\delta \geq 0$;
in other words, $\mathcal{R} (Z) = \{R_\epsilon(Z)\}_{\epsilon}$ is a filtered simplicial complex.
Janko Latschev shows in \cite{latschev2001vietoris} that
when $\X$ is a closed Riemannian manifold,  there is an $\epsilon_0 > 0$, so that  if
$ 0 < \epsilon \leq \epsilon_0$, then there exists $\delta > 0 $ for
which $d_H(X, \X) < \delta$
 implies that
the geometric realization of $R_\epsilon(X)$  is homotopy equivalent to $\X$.
Discarding the manifold hypothesis --- which is not expected to hold in general applications  --- highlights the value of persistence
as a homology inference tool.
Indeed, in \cite{chazal2008towards} Chazal, Oudot and Yan
show that if $\X \subset \R^d$ is compact with positive \emph{weak feature size}\footnote{this is a notion of how complex  the embedding of $\X$ into Euclidean space is.} \cite{chazal2005weak},
 and $X \subset \R^d$ is finite with $d_H(X,\X)$ small,
then there exists a range for $\epsilon > 0$ where
$H^{\epsilon, 3\epsilon}_n(\mathcal{R}(X) ; \F) $ is isomorphic to $ H_n\left(\X^{(\epsilon)}; \F\right)$.
It is worth noting that while these theorems deal with small $\epsilon$,
far less is known about the large-scale regime.
Indeed,  aside from trivial examples,  the circle is  (essentially) the only space $Z$ for which the homotopy type of $R_\epsilon(Z)$
is known explicitly for all $\epsilon > 0 $
 \cite{adamaszek2017vietoris, adamaszek2017vietorisEllipses}.

The efficient computation of the persistent Betti numbers of a finite
filtered  simplicial complex
$\mathcal{K} = \{K_0 \subset K_1 \subset  \cdots \subset K_J = K \}$, was  addressed by
Edelsbrunner, Letscher and Zomorodian  in (2000) \cite{edelsbrunner2000topological},
for   subcomplexes of a triangulated 3-sphere
 and homology with coefficients in  $\F_2 = \{-1,1\}$.
This restriction was a tradeoff between generality and speed:
the algorithm was based on previous work
of Delfinado and Edelsbrunner \cite{delfinado1995incremental}
to compute (standard) Betti numbers incrementally
in time $O(N \alpha(N))$, where $N$ is the number of simplices of $K$ and $\alpha$ is the inverse of the Ackermann function \cite{cormen2009introduction}.
Since the Ackermann function grows very rapidly,
its inverse $\alpha$ grows very slowly.
Though limited in generality, the approach by Delfinado and Edelsbrunner
highlights the following idea:
If  $K_{j}$ is obtained from $ K_{j-1} $ by adding a single simplex $\tau\in K$, and
 $H_n(K_{j-1}; \F) \longrightarrow H_n(K_j;\F)$ is not surjective,
then either $\tau$
is an $n$-simplex creating a new homology class,
or it is an $n+1$-simplex  eliminating a class from $K_{j-1}$.
Thus, simplices in $K$ that either create
or annihilate a given persistent  homology class can be put in pairs $(\tau,\sigma)$ of the form creator-annihilator.
These pairings are in fact a byproduct of the incremental algorithm
of Delfinado and Edelsburnner.
The  \emph{barcode} is also introduced in \cite{edelsbrunner2000topological} as a visualization tool for persistence:
each pair  $(\tau, \sigma)$ yields an interval $[j, \ell)$, where  $j$ (birth time) is the smallest index so that
$\tau \in K_j$, and  $\ell > j$  (death time) is the smallest index  for which  $\sigma \in K_\ell$.
Thus, long intervals indicate stable homological features throughout $\mathcal{K}$,
while short ones reflect  topological noise.
The resulting multiset of intervals (as repetitions are allowed)  is  called a barcode. The notation is
 $\dgm_n(\mathcal{K})$.
Moreover, the barcode subsumes the persistent Betti numbers,  since $\beta_n^{\epsilon, \delta}(\mathcal{K};\F)$
is  the number of intervals $[j,\ell) \in \dgm_n(\mathcal{K})$ with $j \leq \epsilon$ and $\ell > \epsilon +\delta $.
Below in Figure \ref{fig:FilteredComplex} we show an example of a filtered simplicial complex, the simplicial
pairings $(\tau, \sigma)$, and the resulting barcodes.

\begin{figure}[htb!]
  \centering
  \includegraphics[width=\textwidth]{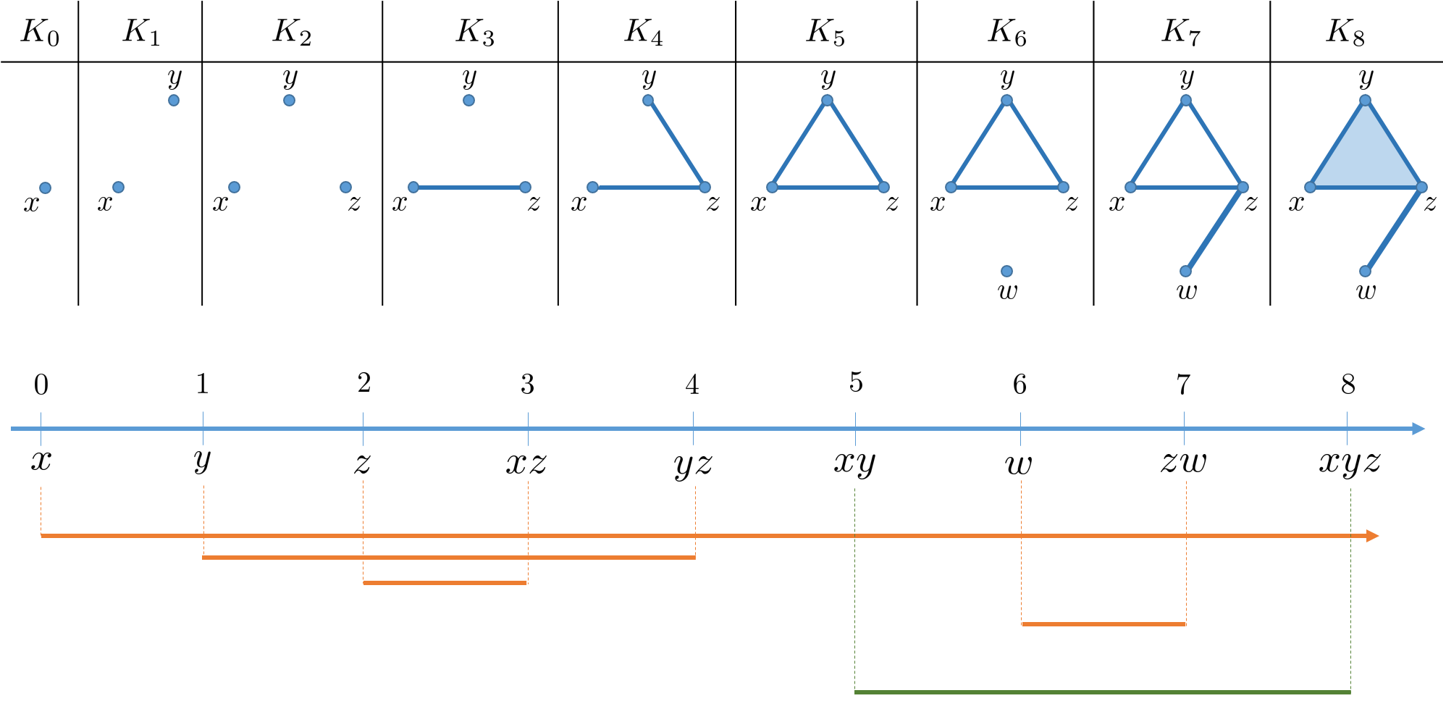}
  \caption{A filtered simplicial complex $\mathcal{K} = \{K_0 \subset \cdots \subset  K_8\}$, along with the simplicial pairings $(\tau, \sigma)$,
  and the resulting barcodes for homology in dimensions 0 (orange) and 1 (green).}
  \label{fig:FilteredComplex}
\end{figure}
\section{Here Comes the  Algebra}

The developments up to this point can be thought of as the computational and geometric era of persistent homology.
 Around 2005 the focus started to shift towards algebra.
Zomorodian and Carlsson introduced in
\cite{zomorodian2005computing}
the persistent homology
\[
PH_n(\mathcal{K}; \F) := \bigoplus_{j \in \Z } H_n(K_j; \F) \;\;\; ,
\;\;\; \mathcal{K} = \{K_j\}_{j \in \Z}
\]
of a filtered complex  $\mathcal{K}$,
as the  graded module over  $\F[t]$
with left multiplication by $t$ on $j$-homogeneous elements
given by the linear map
\[\phi_j : H_n(K_j ; \F) \longrightarrow H_n(K_{j+1};\F)\]
induced by the inclusion $K_j \hookrightarrow K_{j+1}$.
Since then, $PH_n(\mathcal{K};\F)$ is referred to in the literature as a persistence module.
More generally \cite{ bdss:1, bubenik2014categorification},
let $\mathsf{J}$ and $\mathsf{C}$ be categories with $\mathsf{J}$ small  (i.e., so that
its objects  form a set).
The category of $\mathsf{J}$-indexed persistence objects in  $\mathsf{C}$
is defined as the functor category $\mathsf{Fun}(\mathsf{J}, \mathsf{C})$; its objects
are functors $F:\mathsf{J}\rightarrow \mathsf{C}$,
and its morphisms are  natural transformations $\varphi: F_1 \Rightarrow F_2$.
The typical  indexing category comes from having
a partially ordered set $(P,\preceq)$, and letting $\underline{P}$
denote the category with objects $\mathsf{Obj}(\underline{P}) = P$, and  a unique morphism
from $p_1$ to $p_2$ whenever $p_1 \preceq p_2$.
We'll abuse notation and  denote this morphism by $p_1 \preceq p_2$,
instead of the categorical notation $p_1 \rightarrow p_2$.

It can be readily checked that if
$\mathsf{Mod}_R$ denotes the category of (left) modules over a commutative ring $R$ with unity,
and
$g\mathsf{Mod}_{R[t]}$ is  the category of
$\Z$-graded modules over the polynomial ring $R[t]$, then
\begin{equation}\label{eq:PerMod_CatEquiv}
\begin{array}{ccl}
 \mathsf{Fun}(\underline{\Z}, \mathsf{Mod}_R)  &
\longrightarrow &  g\mathsf{Mod}_{R[t]} \\[.2cm]
  M \;\;, \;\; \varphi &  \mapsto & \bigoplus\limits_{j\in \Z}M(j) \;\; , \;\;  \bigoplus\limits_{j\in \Z}\varphi_j
\end{array}
\end{equation}
is an equivalence of categories. On the graded $R[t]$-module side, multiplication by $t$ on $j$-homogeneous elements
is given by $M( j \leq j+1) : M(j) \longrightarrow M(j+1)$.
This equivalence  shows why/how the evolution of homological features in a $\Z$-filtered complex $\mathcal{K}$,
can be encoded as the algebraic structure of the persistence module $PH_n(\mathcal{K};\F)$.

\subsection{Persistence Modules and Barcodes}
When $PH_n(\mathcal{K};\F)$ is finitely generated as an $\F[t]$-module --- e.g.
 if  $K_j = \emptyset$ for $j < 0$ and $\bigcup\limits_{j \in \Z} K_j$ is finite ---
then one has a graded isomorphism
\begin{equation}\label{eq:PerMod_Decomp}
PH_n(\mathcal{K};\F) \cong
\left(
\bigoplus_{q=1}^{Q}
t^{n_q} \cdot
\F[t]
\right)
\oplus
\left(
\bigoplus_{\ell=1}^{L}
\left(t^{m_\ell} \cdot
\F[t]\right)/(t^{  m_\ell + d_\ell})
\right)
\end{equation}
for $n_q, m_\ell \in  \Z$ and $d_\ell \in \N $  \cite{webb1985decomposition}.
The  decomposition (\ref{eq:PerMod_Decomp}) is unique up to permutations, and thus the
 intervals
\begin{align*}
[n_1,\infty),[n_2,\infty), \ldots, & [n_Q,\infty), \\ &[m_1, m_1 + d_1), [m_2,m_2 + d_2)
,\ldots, [m_L , m_L + d_L)
\end{align*}
provide  a complete discrete invariant for (i.e., they uniquely determine)  the $\F[t]$-isomorphism type of $PH_n(\mathcal{K};\F)$.
Moreover, this multiset
recovers the barcode $\dgm_n(\mathcal{K})$ of Edelsbrunner, Letscher and Zomorodian \cite{edelsbrunner2000topological}.%

Carlsson and Zomorodian also observe that $PH_n(\mathcal{K};\F)$ is in fact the  homology of an appropriate  chain complex of graded $\F[t]$-modules.
Hence, a graded version of the
 Smith Normal Form  \cite{dumas2003computing}  computes
the barcode decomposition (\ref{eq:PerMod_Decomp}), providing a general-purpose algorithm.
This opened the flood gates;
barcodes could now be computed as a linear algebra problem   for any
finite filtered simplicial complex
$ K_0 \subset \cdots \subset K_J = K $, over any (in practice finite) field
of coefficients, and up to any homological dimension.
The resulting matrix reduction algorithm,  implemented initially in the  JPlex library (now javaPlex) \cite{adams2014javaplex},
runs in polynomial time: its worst time complexity is
$O(N^3)$, where $N$ is the number of simplices of $K$.
In fact Dmitriy Morozov exhibits in \cite{morozov2005persistence}  a finite filtered complex of dimension 2, attaining the worst-case.
This shows  that the cubic bound is tight for general barcode computations.

While this sounds potentially slow, specially compared to the time complexity $O(N\cdot \alpha(N))$ of the sequential algorithm,
Morozov's example should be contrasted with filtrations arising from applications.
In practice the matrices to be reduced are sparse, and
computing their associated barcode decomposition takes at worst matrix-multiplication time
$O(N
^{\omega})$
\cite{milosavljevic2011zigzag}, where  $\omega \approx 2.373 $ \cite{williams2012multiplying}.
Over the last ten years or so there has been a flurry of activity towards better implementations  and
faster persistent homology computations.
A recent survey \cite{otter2017roadmap} compares several leading open source libraries for computing persistent homology.
All of them implement different optimizations,  exploit new theoretical developments and novel heuristics/approximations.
For instance, one improvement is to first simplify the input filtered complex without changing its persistent homology
(e.g., using discrete Morse theory \cite{mischaikow2013morse});
or to compute persistent cohomology, since it is more efficient than persistent homology and gives the same answer \cite{de2011dualities}.

The shift towards algebra has had other important consequences; specifically:
1) a better understanding of stability for barcodes,  and
2) several theorems  describing how the choice of  categories $\mathsf{J} $ and $\mathsf{C}$ impacts the  computability of  isomorphism invariants for objects in $\mathsf{Fun}(\mathsf{J}, \mathsf{C})$.
Let me say a few words about stability.

\subsection{The Stability of Persistence}
Let $\X$ be a triangulable topological space (e.g., a smooth manifold)
and let $f: \X \longrightarrow \R$ be a tame function (this is a generalization of being  Morse).
The prototypical example in TDA arises from a compact set  $X\subset \R^d$, and letting $f_X : \R^d \longrightarrow \R$ be
 $$f_X(y) = \inf\limits_{x\in X} \|x - y\|.$$
 Thus $f_X^{-1}(-\infty, \epsilon] = X^{(\epsilon)}$.
Given $f: \X \longrightarrow \R$, let $\dgm_n(f)$ denote the barcode for the $n$-th persistent homology
of $\left\{f^{-1}(-\infty,\epsilon]\right\}_{\epsilon\in \R}$.
Drawing inspiration from Morse theory, Cohen-Steiner, Edelsbrunner and
Harer introduced in  (2007) \cite{cohen2007stability}  two foundational ideas:
 (1)
the bottleneck distance $d_B(\,\cdot\, ,\,\cdot\,)$ between  barcodes,
and (2) the stability theorem   asserting that for tame $f,g : \X \longrightarrow \R$ one has that\footnote{A similar result was established in \cite{d2003optimal} for $n=0$.}
\[
d_B(\dgm_n(f), \dgm_n(g)) \leq \|f - g\|_\infty
\]
In particular,  if $X, Y \subset \R^d$
are compact and $\mathcal{X} = \left\{X^{(\epsilon)}\right\}_\epsilon$, $\mathcal{Y}= \left\{Y^{(\epsilon)}\right\}_\epsilon$,
then $d_B(\dgm_n(\mathcal{X}), \dgm_n(\mathcal{Y})) \leq d_H(X, Y)$.
This inequality implies that slight changes  to the input data  change  the barcodes slightly,
which is key for applications where (Hausdorff) noise plays a role.

Towards  the end of 2008 Chazal et.~al.~solidified the idea of stability with the introduction of
interleavings for $\R$-indexed persistence modules \cite{chazal2009proximity}.
The construction is as follows.
For $\delta \geq 0$
let $T_\delta : \underline{\R} \longrightarrow  \underline{\R} $
be the translation functor $T_\delta(\epsilon) = \epsilon + \delta$.
An $\delta$-interleaving between
two  persistence vector spaces
$V,W:\underline{\R}\longrightarrow \mathsf{Mod}_\F$
is a pair $(\varphi,\psi)$ of natural transformations
\[
\varphi: V \Rightarrow W\circ T_{\delta}
\;\;\;\;\;\;\; \mbox{ and } \;\;\;\;\;\;\;
\psi: W \Rightarrow V\circ T_\delta
\]
so that
$\psi_{\epsilon + \delta} \circ \varphi_\epsilon = V(\epsilon \leq \epsilon + 2\delta)$ and
$\varphi_{\epsilon + \delta} \circ \psi_\epsilon = W(\epsilon \leq \epsilon + 2\delta)$ for all $\epsilon \in \R$.
The interleaving distance between $V$ and $W$, denoted $d_I(V,W)$, is  defined
as the infimum over all $\delta \geq 0$ so that $V$ and $W$ are $\delta$-interleaved,
if interleavings  exist.
If there are no interleavings, the distance is defined as
$\infty$.
It readily follows that $d_I$ is an extended (since infinity can be a value) pseudometric on $\mathsf{Fun}(\underline{\R}, \mathsf{Mod}_\F)$, and that $d_I(V,W) = 0$ whenever $V\cong W$.
The converse, however,
is false in general (more on this later).

Chazal et. al. \cite{chazal2009proximity} show that if $V: \underline{\R} \longrightarrow \mathsf{Mod}_\F$  satisfies $\mathsf{rank}\big( V(\epsilon < \epsilon')\big)  < \infty$ for all pairs $\epsilon < \epsilon'$,
this is called being $\mathsf{q}$-tame,
then $V$ has a well-defined barcode $\dgm(V) $
(see \cite{crawley2015decomposition} for a shorter proof when $\mathsf{dim}_\F V(\epsilon) < \infty$ for all $\epsilon$; this is called being pointwise finite).
Moreover,  if $V,W$ are $\mathsf{q}$-tame, then one has the algebraic stability theorem
$d_B(\dgm(V), \dgm(W)) \leq d_I(V,W)$.
This  turns out to be an equality:
\[
d_B(\dgm(V), \dgm(W)) = d_I(V,W)
\]
which nowadays is referred to as the Isometry Theorem;
the first known proof is due to  Lesnick \cite{lesnick2015theory}.

As  I said  earlier,  $d_I(V,W)$ can be zero  for $V $  and $W$  nonisomorphic,
and thus  $\dgm(V)$ is not a complete  invariant
in the $\mathsf{q}$-tame $\R$-indexed case.
This can be remedied as follows.
Let $\mathsf{qFun}(\underline{\R}, \mathsf{Mod}_\F)$
denote the full subcategory of $\mathsf{Fun}(\underline{\R},\mathsf{Mod}_\F)$
comprised of
 $\mathsf{q}$-tame persistence modules.
The ephemeral category  $\mathsf{eFun}(\underline{\R}, \mathsf{Mod}_\F)$,
is the full subcategory of $\mathsf{qFun}(\underline{\R},\mathsf{Mod}_\F)$
with objects $V : \underline{\R} \longrightarrow \mathsf{Mod}_\F$
satisfying $V(\epsilon < \epsilon' ) = 0$ for all
$\epsilon < \epsilon'$.
The
observable category
$\mathsf{oFun(\underline{\R},\mathsf{Mod}_\F)}$
is the quotient category
\[
\mathsf{qFun}(\underline{\R},\mathsf{Mod}_\F)/ \mathsf{eFun}(\underline{\R},\mathsf{Mod}_\F)
\]
As shown by Chazal et. al. in \cite{chazal2014observable},
$d_I$ descends to an extended metric on the observable category,
and taking barcodes
\[
\dgm : \big(\mathsf{oFun}(\underline{\R}, \mathsf{Mod}_\F) , d_I\big) \longrightarrow \big(\Dgm , d_B\big)
\]
is an isometry.
Hence, the barcode is a  complete invariant for the
isomorphism type of observable $\R$-indexed persistence
vector spaces.
In summary, the modern view of stability is algebraic; persistence modules are compared via
interleavings, which one then tries to relate to the bottleneck distance between the associated barcodes if they exist.

\subsection{Changing Indexing Categories: Multi-d Persistence, Quivers and Zigzags}
One of the early realizations in TDA
was the usefulness  of having filtrations indexed by more than one parameter (1999) \cite{frosini1999size}.
For instance, given a  data set $X\subset \R^d$ one might want
to focus on   densely-populated regions \cite{carlsson2008local},
or portions with high/low curvature \cite{carlsson2005persistence}.
This leads naturally to $\Z^n$-filtered complexes: $\{K_\uu\}_{\uu \in \Z^n }$,
$\uu = (u_1,\ldots, u_n)$, where $K_{\uu} \subset K_{\vv}$ whenever
$\uu \preceq \vv$ (i.e.,
$u_1 \leq v_1, \ldots, u_n \leq v_n$).
In this multi-filtered complex,
each  filtering direction $u_1,\ldots, u_n$
is meant to capture an attribute: e.g. distance/scale, density, curvature, etc.
Taking homology with coefficients in $\F$ yields objects in
 $\mathsf{Fun}(\underline{\Z^n}, \mathsf{Mod}_\F)$,
 and just like before,
$\Z^n$-indexed persistence $\F$-vector spaces
correspond to   $n$-graded modules over the $n$-graded polynomial ring
$P_n := \F[t_1,\ldots, t_n]$.
Parameterizing the  isomorphism classes of said modules,
 for $n\geq2$, turns out to be much more involved
than the barcodes from $n=1$.
Indeed, around 2009 Carlsson and Zomorodian  \cite{carlsson2009theory} showed
that the isomorphism type of a finitely generated $n$-graded $P_n$-module
is uniquely determined by the following data:
two finite multisets $\xi_0,\xi_1 \subset \Z^n$
encoding the location and multiplicity of birth-death events, and a point in the quotient of an algebraic variety
$\mathcal{RF}(\xi_0,\xi_1)$  by
the algebraic action of an algebraic group $G_{\xi_0}$.
The multisets $\xi_0, \xi_1$
are the discrete portions of the resulting isomorphism invariant, while
$\mathcal{RF}(\xi_0,\xi_1)/ G_{\xi_0} $ parameterizes the (potentially) continuous part.
Here is an example due to Carlsson and Zomorodian \cite{carlsson2009theory}
illustrating how complicated this quotient can be.
For $n=2$, consider the isomorphism classes of $P_2$-modules having
$\xi_0 = \{(0,0), (0,0)\}$ and
$\xi_1 = \{(3,0),(2,1), (1,2), (0,3)\}$.
If
$\mathsf{Gr}_1(\F^2)$ denotes the Grassmannian of lines
in $\F^2$, then
\[
\mathcal{RF}(\xi_0,\xi_1)  =  \mathsf{Gr}_1(\F^2)
\times
\mathsf{Gr}_1(\F^2)
\times
\mathsf{Gr}_1(\F^2)
\times
\mathsf{Gr}_1(\F^2)
\] and
$G_{\xi_0} $ turns out to be the degree 2 general linear group  $\mathsf{GL}_2(\F)$
acting diagonally on $\mathsf{Gr}_1(\F^2)^4$.
Since $\mathsf{Gr}_1(\F^2)^4 / \mathsf{GL}_2(\F)$  contains a copy of
$\F \smallsetminus \{0,1\}$,
and each point in this set yields a distinct isomorphism class of
 $P_2$-modules,
it follows that there is no
complete discrete invariant for (finite!) multi-d persistence.

The vast majority of recent results from multidimensional persistence
focus on  computable descriptors/visualizations
of its intricate algebraic structure.
 Besides introducing   the  parametrization $$\xi_0,\xi_1, \mathcal{RF}(\xi_0,\xi_1)/G_{\xi_0},$$
Carlsson and Zomorodian also propose the rank invariant:
For a $\mathsf{q}$-tame module $V :\underline{\Z^n} \longrightarrow \mathsf{Mod}_\F$, it is defined as the function
 $\rho_V$   sending each pair  $\uu \preceq \vv  $
in $\Z^n$
to the integer $\mathsf{rank} \, V(\uu \preceq \vv)$.
$\rho_V$ is  computable (see \cite{carlsson2009computing} for a polynomial-time algorithm),
it is discrete, and an invariant of the isomorphism
type of $V$.
When $n=1$ one can recover $\dgm(V)$  from
$\rho_V$ and viceversa,
and thus $\rho_V$ is complete in the 1-dimensional case.
Knudson notes in \cite{knudson2008refinement}
that $\xi_0(V)$ and $\xi_1(V)$ are in fact the locations/multiplicities
of birth events  in the torsion modules
$\mathsf{Tor}^{P_n}_0 (V,\F)$ and $\mathsf{Tor}^{P_n}_1(V,\F)$,
respectively; here
$\F$ is identified with the $P_n$-module 
\[
\F[t_1,\ldots, t_n]/(t_1,\ldots, t_n)\]
The higher-dimensional analogs
$\mathsf{Tor}^{P_n}_j(V,\F)$, $j\geq 2$, lead to
a family of finite multisets $\xi_j(V) \subset \Z^n$, each with its
own geometric interpretation,
serving as isomorphism invariants for $V$.
Other approaches to invariants for multidimensional persistence include  the Hilbert Series of
Harrington et. al.  \cite{harrington2017stratifying},
the extended algebraic functions of
Skryzalin and Carlsson \cite{skryzalin2017numeric},
and the feature counting invariant of
Scolamiero et. al. \cite{scolamiero2017multidimensional}.
Lesnick and Wright have recently released
RIVET, the Rank Invariant Visualization and Exploration Tool \cite{lesnick2015interactive}.
Put simply, RIVET uses the fact that
 if $ V : \underline{\R^2} \longrightarrow \mathsf{Mod}_\F$
is $\mathsf{q}$-tame and  $L\subset \R^2$ is a line with nonnegative slope (hence
a totally ordered subset of $(\R^2,\preceq)$), then  $V^L : \underline{L} \longrightarrow \mathsf{Mod}_\F$,  the 1-dimensional  persistence vector space
obtained  by restricting
$V$ to $L$, has a well-defined barcode $\mathsf{bc}\left(V^L\right)$.
The key feature in RIVET is a graphical interface which, for finite bi-filtrations, displays
 $\mathsf{bc}\left(V^L\right)$ interactively as the user varies  $L$.
This is particularly useful for parameter selection and the exploratory
analysis of  data sets with filtering functions.

Multidimensional persistence is a great example of how
a seemingly innocuous change in indexing category, say from
$\underline{\Z}$ to $\underline{\Z^2}$, can lead to a widely different and much more complicated classification problem.
With this in mind, one would like to have a systematic approach to
address  the ensuing complexity.
The representation theory of Quivers \cite{derksen2005quiver} offers one such
avenue.
It turns out that the  classification of finite
$\mathsf{J}$-indexed persistence vector spaces $ V: \mathsf{J} \longrightarrow \mathsf{Mod}_\F$
can be studied directly from the shape of the indexing category $\mathsf{J}$.
Indeed, let $G(\mathsf{J})$ be the finite directed (multi)graph with the objects of  $\mathsf{J}$ as vertices, and one arrow for
every morphism that is neither an identity nor a  composition.
Also, let $\gor{G}(\mathsf{J})$ be the undirected graph obtained
from $G(\mathsf{J})$ by forgetting arrow directions.
When $\gor{G}(\mathsf{J})$ is acyclic,
Gabriel's theorem \cite{gabriel1972unzerlegbare} implies
that pointwise finite objects in
$\mathsf{Fun}(\mathsf{J}, \mathsf{Mod}_\F)$ can be classified via  complete discrete invariants,
 if and only if  the connected components of $\gor{G}(\mathsf{J})$
 are Dynkin diagrams of the types described in Figure \ref{fig:Dinkyn} below.

\begin{figure}[htb!]
  \centering
  \includegraphics[width=0.7\textwidth]{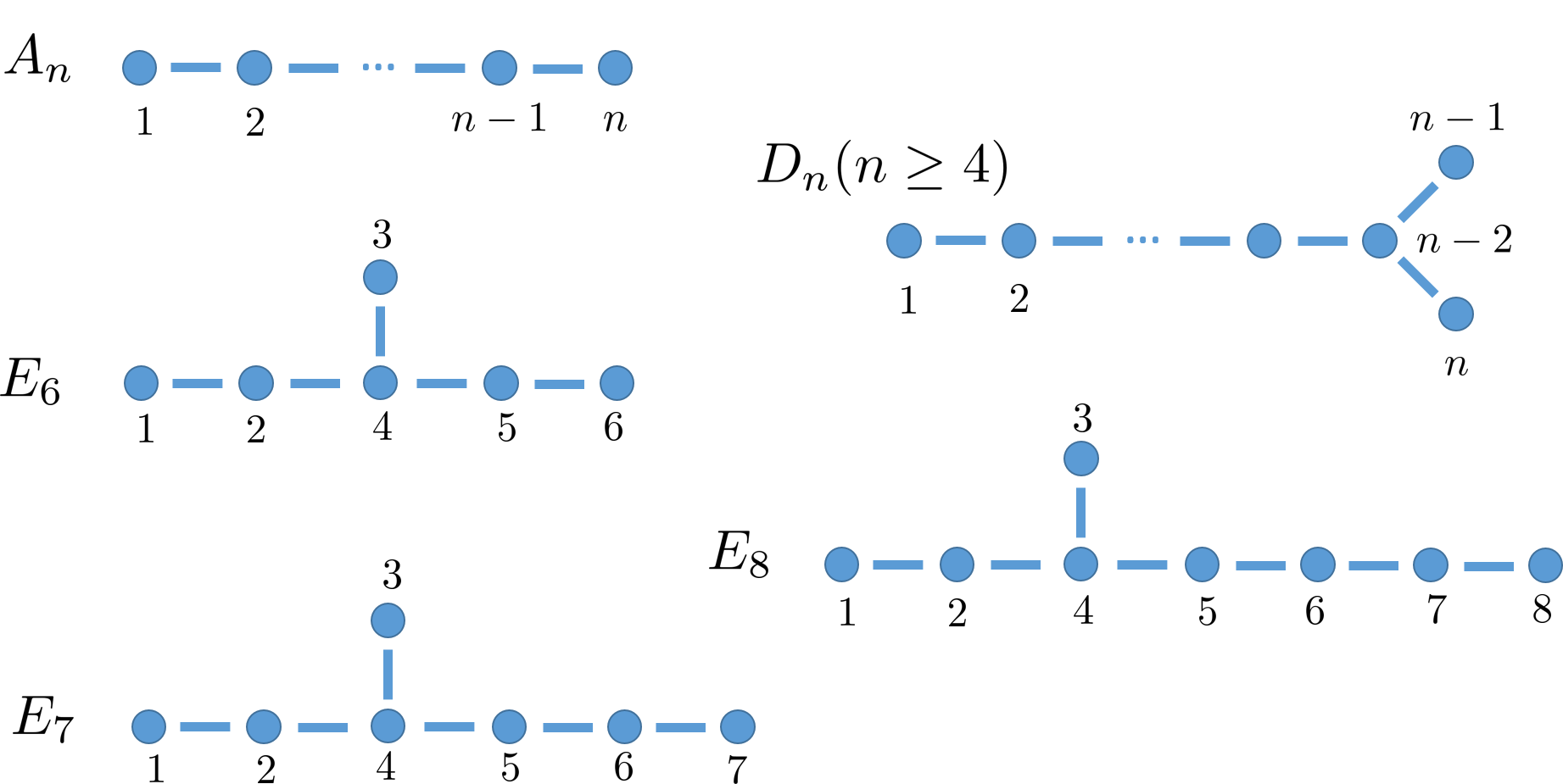}
  \caption{Dynkin diagrams of type
$A_n$ for $n\geq 1$, $D_n$ for $n\geq 4$, and $E_n$ for $n = 6,7,8$.}
\label{fig:Dinkyn}
\end{figure}

Here is an example of how this result can be used to avoid unpleasant surprises: Suppose that $G(\mathsf{J})$ is the graph with vertices $x_0 ,\ldots ,x_N$
and $N\geq 5$  edges $x_n\rightarrow x_0$, $n =1,\ldots, N$  (see Examples 3 and 8 in \cite{derksen2005quiver}).
While the resulting $\mathsf{J}$-indexed persistence vector spaces
$V : \mathsf{J} \longrightarrow \mathsf{Mod}_\F$ may look simple (just star-shaped, right?),
the connected graph $\gor{G}(\mathsf{J})$ is not a Dynkin diagram,
and the ensuing classification problem is in fact of ``wild type'':
complete invariants must include continuous high-dimensional pieces, just like in  multidimensional persistence.

These ideas entered the TDA lexicon around 2010
 with the definition of Zigzag persistence by Carlsson and de Silva  \cite{carlsson2010zigzag}.
Regular persistence addresses the problem of identifying
stable homological features in a monotone system of spaces and continuous maps $$X_1 \rightarrow X_2 \rightarrow \cdots \rightarrow X_J.$$
Zigzag persistence, on the other hand, is a generalization to the non-monotone case.
Here is a practical example: suppose one has an ordered sequence
of spaces $X_1,\ldots, X_J$ (e.g., from time varying data), but no obvious maps $X_j \rightarrow X_{j+1}$.
The need to track topological features as $j$ varies leads
one   to consider the system
\[
X_1 \hookrightarrow X_1 \cup X_2 \hookleftarrow X_2 \hookrightarrow \cdots \hookleftarrow X_j \cup X_{j+1} \hookrightarrow \cdots \hookleftarrow X_J
\]
and the resulting \emph{zigzag diagram}
$$V_1 \rightarrow V_2 \leftarrow V_3 \rightarrow \cdots \leftarrow V_n$$
 at the homology level.
More generally, a (finite) zigzag is a sequence of vector spaces
$V_1,\ldots, V_n$ and linear maps $V_j \rightarrow V_{j+1}$ or
$V_{j} \leftarrow V_{j+1}$. The sequence of arrow directions,
e.g. $\tau = (\mathsf{left,left,right, . . . , right,left})$, is  the  zigzag type.
Since in this case
any  choice of $\tau$ forces the
indexing category $\mathsf{J}_\tau$  to satisfy
$\gor{G}(\mathsf{J}_\tau) = A_n$ (one of the aforementioned Dynkin diagrams),
then Gabriel's theorem implies
that finite zigzags $$V: \mathsf{J}_\tau \longrightarrow \mathsf{Mod}_\F$$
are completely classified by a discrete invariant.
Just as for regular 1-dimensional persistence the invariant turns out to be a barcode,
which can be efficiently computed \cite{milosavljevic2011zigzag},
and for which there is  a zigzag stability theorem \cite{botnan2016algebraic}
recently established by Botnan and Lesnick.

When the graph $\gor{G}(\mathsf{J})$ has cycles,
the functoriality of objects in $\mathsf{Fun}(\mathsf{J}, \mathsf{Mod}_\F)$
is captured by the notion of a quiver with relations.
The taxonomy from Gabriel's theorem no longer applies, but one can still find some answers
in the representation theory of associative algebras.
A particularly important instance  is
 when the cycles of $\gor{G}(\mathsf{J})$ are not oriented cycles in $G(\mathsf{J})$;
 in this case the  algebras of interest are finite dimensional (hence Artinian)
and
Auslander-Rieten theory \cite{auslander1997representation} becomes relevant.
Escolar and Hiraoka \cite{escolar2016persistence} have recently put these
ideas to use in the context of persistent objects indexed by commutative ladders; that is, the persistence of a morphism between two zigzags
of the same type:
\[
\begin{tikzcd}
\bullet \ar{r} \ar{d}& \bullet \ar{d}& \bullet \ar{l}\ar{d} & \ar{l}\ar{r}\cdots &  \bullet \ar{r}\ar{d}& \bullet \ar{d}& \bullet \ar{d} \ar{l}\\
\bullet \ar{r} & \bullet & \bullet \ar{l} & \ar{l}\ar{r}\cdots &  \bullet \ar{r} & \bullet & \bullet \ar{l}
\end{tikzcd}
\]

The resulting theory sits somewhere between zigzag persistence and multi-dimen\-sional persistence: short ladders (length $\leq $ 4) have complete discrete invariants,
but longer ones do not.
Escolar and Hiraoka present an algorithm for computing these invariants, and also an interesting application to computational chemistry.

I think this is a good place for me to stop; hopefully it is also a good starting point for the reader
interested in persistent homology.
There are several books covering many of the ideas I presented here,
as well as many others.
The interested reader would certainly benefit from these resource \cite{ghrist2014elementary, chazal2016structure, edelsbrunner2010computational, oudot2015persistence}.

\end{document}